\documentclass[a4paper,11pt]{amsart}

\usepackage{amsmath,amssymb,amsfonts,amsthm}

\usepackage[alphabetic]{amsrefs}
\usepackage{bm}
\usepackage{graphicx}
\usepackage{ascmac}
\usepackage[all]{xy}
\usepackage{amsthm}
\usepackage{tikz}
\usetikzlibrary{patterns}
\usepackage{pxpgfmark}
\usepackage{multirow}
\usetikzlibrary{intersections,arrows,calc,decorations.markings} 
\newtheorem{theorem}{Theorem}[section]
\newtheorem{proposition}[theorem]{Proposition}

\newtheorem{corollary}[theorem]{Corollary}
\newtheorem{lemma}[theorem]{Lemma}

\theoremstyle{definition}
\newtheorem{remark}[theorem]{Remark}

\newtheorem{definition}[theorem]{Definition}

\makeatletter
 
 \@addtoreset{equation}{section}
\makeatother

\newcommand{\opname}[1]{\operatorname{\mathsf{#1}}}
\newcommand{\Hom}{\opname{Hom}}

\newcommand{\Ext}{\opname{Ext}}
\newcommand{\dimv}{\opname{\underline{dim}}}

\title[homological properties of $\opname{rep}(Q,\mathbb{F}_1)$]{On homological properties of the category of $\mathbb{F}_1$-representations over a linear quiver of type $\mathbb{A}_n$}\thanks{Partially supported by the National Natural Science Foundation of China (Grant No. 11971326)}
\author{Changjian Fu}
\author{Longjun Ran}
\author{Liang Yang}
\address{Department of Mathematics, Sichuan University, Chengdu, 610064 PR China}

\email{changjianfu@scu.edu.cn(Fu)}
\email{1197882950@qq.com(Ran)}
\email{malyang@scu.edu.cn(Yang)}

\dedicatory{Dedicated to Professor Liangang Peng on the Occasion of his 65th Birthday}
\date{}
\begin{document}

\maketitle

\begin{abstract}
    Let $Q$ be a quiver of type $\mathbb{A}_n$ with linear orientation and $\opname{rep}(Q,\mathbb{F}_1)$ the category of representations of $Q$ over the virtual field $\mathbb{F}_1$.
It is proved that $\opname{rep}(Q,\mathbb{F}_1)$ has global dimension $2$ whenever $n\geq 3$ and it is hereditary if $n\leq 2$.  As a consequence, the Euler form $\langle L, M\rangle=\sum_{i=0}^\infty (-1)^i\dim \Ext^i(L,M)$ is well-defined. However, it does not descend to the Grothendieck group of $\opname{rep}(Q,\mathbb{F}_1)$. This yields  negative answers to  questions raised by Szczesny in [IMRN, Vol. 2012, No. 10, pp. 2377–2404].
\end{abstract}

\tableofcontents

\section{Introduction}
Let $\mathbb{F}_1$ be the virtual field, i.e., the field with one element. Quiver representations over $\mathbb{F}_1$ were introduced by Szczesny \cite{Sz12} as a degenerated combinatorial model of quiver representations over a field. It has been shown by Szczesny \cite{Sz12} that the category $\opname{rep}(Q,\mathbb{F}_1)$ of representations of a quiver $Q$ over $\mathbb{F}_1$ has  many nice properties like the one $\opname{rep}(Q,k)$, the category of representations of $Q$ over a field $k$. In particular, it has a zero object, kernels, co-kernels, and direct sums. The
Jordan–H\"{o}lder  and Krull–Schmidt theorems also hold for $\opname{rep}(Q,\mathbb{F}_1)$.

An evident difference between $\opname{rep}(Q,\mathbb{F}_1)$ and $\opname{rep}(Q,k)$ is that the category $\opname{rep}(Q,\mathbb{F}_1)$ is not additive. 
 This leads to tremendous differences between $\opname{rep}(Q,\mathbb{F}_1)$ and $\opname{rep}(Q,k)$. For instance, it is well-known that a connected quiver $Q$ is of representation-finite over a field if and only if the underlying diagram of $Q$ is of type $\mathbb{ADE}$, see \cite{Ga72}. However, it was proved by Jun and Sistko \cite{JS23} that $Q$ is of representation-finite over $\mathbb{F}_1$ if and only if $Q$ is a tree quiver. Nevertheless, the category $\opname{rep}(Q,\mathbb{F}_1)$ is a proto-exact category in the sense \cite{DK12}, which is a non-additive generalization of Quillen's exact category. Hence the Hall algebra theory of $\opname{rep}(Q,\mathbb{F}_1)$ can be established as usual.
In \cite{Sz12}, Szczesny has investigated Hall algebras of $\opname{rep}(Q,\mathbb{F}_1)$ for certain classes of finite quivers. In particular, the nilpotent halves  of  Lie algebras  of finite type $\mathbb{A}$ and affine type $\hat{\mathbb{A}}$ and their enveloping algebras can be realized as Hall algebras of $\opname{rep}(Q,\mathbb{F}_1)$ for  quivers of type $\mathbb{A}$ and cyclic quivers  respectively. Such a realization can be understood as a degeneration of Ringel's realization \cite{R90} of the positive part of a quantum group  as Hall algebra for quiver representations over a finite field. 
 The Euler form \[\langle L,N\rangle:=\sum_{i=0}^\infty(-1)^i\dim_k\Ext^i(L,N),\] 
 where $L,N\in \opname{rep}(Q,k)$, has played an important role in Ringel's realization and also in the study of quiver representations over $k$. It is expected by Szczesny that the Euler form may also play a similar role in the study of representations over $\mathbb{F}_1$. In particular, the following questions should be considered first, see \cite[Section 12]{Sz12}.
 \begin{itemize}
     \item[(a)] Is the category $\opname{rep}(Q,\mathbb{F}_1)$ hereditary?
     \item[(b)] Is the Euler form well-defined?
     \item[(c)] If the Euler form is well-defined, does it descend to the Grothendieck group of $\opname{rep}(Q,\mathbb{F}_1)$?
 \end{itemize}
Since $\opname{rep}(Q,\mathbb{F}_1)$ is not additive,  traditional homological tools can not be applied. It seems likely that it is very hard to obtain answers for these questions for a general quiver. 

This note is devoted to answering the above questions for a special class of quivers: the linear quiver of type $\mathbb{A}_n$. It turns out that $\opname{rep}(Q,\mathbb{F}_1)$ is not hereditary in general, which reveals new difference between $\opname{rep}(Q,\mathbb{F}_1)$ and $\opname{rep}(Q,k)$. Furthermore, even the Euler form is well-defined, it may not descend to the Grothendieck group of $\opname{rep}(Q,\mathbb{F}_1)$.

The paper is organized as follows. In Section \ref{s:basic-rep}, we recall basic definitions of $\mathbb{F}_1$-representation over a finite quiver and the Yoneda's extension. In Section \ref{s:linear-quiver}, we present the main results of this note. In particular, the global dimension of $\opname{rep}(Q,\mathbb{F}_1)$ for a linear quiver $Q$ of type $\mathbb{A}_n$ is obtained. Furthermore, we also obtain a characterization of projective representations of $\opname{rep}(Q,\mathbb{F}_1)$. In Section \ref{s:exam}, by computing an example, we point out that the Euler form does not descend to the Grothendieck group. We also present an example to illustrate that  part of results in Section \ref{s:linear-quiver} can not be generalized to an arbitrary quiver.

\noindent{\bf Acknowledgments.} We are grateful to Markus Kleinau for pointing out the example in Section \ref{s:non-linear} to indicate that the direct sum of projectives is not projective in general. The authors thank Haicheng Zhang and the anonymous referees for suggesting improvements
to the exposition.

\section{Preliminaries}\label{s:basic-rep}
In this section, we recall basic definitions and properties concerning representations of quivers over the virtual field $\mathbb{F}_1$.  We follow \cites{Sz12, JS23} and refer to \cite{Sz12} for unexplained definitions related to representations of quivers over $\mathbb{F}_1$.

\subsection{Vector space over $\mathbb{F}_1$ }

A finite dimensional {\it $\mathbb{F}_1$-vector space} is a finite pointed set $V:=(V,0_V)$. We say that $V$ is of dimension $\dim V:=|V|-1$. An {\it $\mathbb{F}_1$-linear map} from $V:=(V,0_V)$ to $W:=(W,0_W)$  is a pointed function $f:V\to W$  such that $f|_{V\backslash f^{-1}(0_W)}$ is an injection. 
We denote by $f^t: W\to V$ the {\it duality} of $f$, which is the pointed function defined by $f^t(w)=f^{-1}(w)$ for $0_W\neq w\in \opname{im} f=f(V)$ and $f^t(w)=0_V$ otherwise.

Denote by $\text{Vect}(\mathbb{F}_1)$ the category of finite dimensional vector spaces over $\mathbb{F}_1$, whose objects are finite dimensional vector spaces over $\mathbb{F}_1$, while morphisms are $\mathbb{F}_1$-linear maps. The category $\text{Vect}(\mathbb{F}_1)$ has almost all the good properties of the category of vector spaces over a field except that the homomorphism space $\Hom(V,W)$ from $V$ to $W$ has no additive structure, which is only a pointed set, equivalently, an $\mathbb{F}_1$-vector space. We refer to \cites{Sz12,JS23} for details. Following \cite{Sz12}, we denote by $\mathbf{k}$ the one-dimensional $\mathbb{F}_1$-vector space.


\subsection{Representations of quivers over $\mathbb{F}_1$}
A quiver is a quadruple $Q=(Q_0,Q_1,s,t)$ consisting of 
\begin{itemize}
    \item a set of vertices $Q_0$;
    \item a set of arrows $Q_1$;
    \item two maps $s,t:Q_1\to Q_0$, which assign to each arrow $\alpha\in Q_1$ its source $s(\alpha)$ and target $t(\alpha)$.
\end{itemize}
The quiver $Q$ is finite if both $Q_0$ and $Q_1$ are finite sets.

Let $Q$ be a finite quiver. A representation over $\mathbb{F}_1$ is a collection $M=(M_i, M_\alpha)_{i\in Q_0,\alpha\in Q_1}$, where
\begin{itemize}
    \item[(1)] $M_i$ is an $\mathbb{F}_1$-vector space for each vertex $i\in Q_0$;
    \item[(2)] $M_\alpha$ is an $\mathbb{F}_1$-linear map for each arrow $\alpha\in Q_1$.
\end{itemize}
The integer vector \[\dimv M=(\dim M_i)_{i\in Q_0}\in \mathbb{Z}^{Q_0}\] is called the dimension vector of the representation $M$.

For representations $M=(M_i,M_\alpha)$ and $N=(N_i,N_\alpha)$ of $Q$ over $\mathbb{F}_1$, a morphism $\Phi:M\to N$ is a collection of $\mathbb{F}_1$-linear maps $\Phi=(\phi_i)_{i\in Q_0}$, where $\phi_i:M_i\to N_i$, such that the following diagram commutes for each arrow $\alpha:i\to j\in Q_1$:
\[
\xymatrix{M_i\ar[r]^{M_\alpha}\ar[d]^{\phi_i}&M_j\ar[d]^{\phi_j}\\
 N_i\ar[r]^{N_\alpha}&N_j.}
\]
The kernel $\ker \Phi$ and image $\opname{im} \Phi$ of $\Phi$ can be defined as usual, which are subrepresentations of $M$ and $N$ respectively (cf. \cite[Section 4]{Sz12}).

Denote by $\opname{rep}(Q,\mathbb{F}_1)$ the category of representations of $Q$ over $\mathbb{F}_1$. The category $\opname{rep}(Q,\mathbb{F}_1)$ shares many good properties of the category $\opname{rep}(Q,k)$, i.e., the category of representations of $Q$ over a field $k$. In particular, it has a zero object, kernels, co-kernels, and direct sums, which inherits from the category $\text{Vect}(\mathbb{F}_1)$. However, the homomorphism space $\Hom(M,N)$ in $\opname{rep}(Q,\mathbb{F}_1)$ has no additive structure as $\text{Vect}(\mathbb{F}_1)$. Nevertheless, Szczesny \cite{Sz12} proved  that the
Jordan–H\"{o}lder and Krull–Schmidt theorems hold in $\opname{rep}(Q,\mathbb{F}_1)$. Moreover, the classification of indecomposable representations over a tree quiver is obtained.
\begin{theorem}\cite[Theorem 5.1]{Sz12}\label{t:class-tree}
    Let $Q$ be a tree quiver. There is a one-to-one correspondence between the set of connected subquivers of $Q$ and the set of isomorphism classes of indecomposable representations of $Q$ over $\mathbb{F}_1$. Furthermore, for  a connected subquiver $Q'$ of $Q$, the indecomposable representation $M=(M_i, M_\alpha)_{i\in Q_0, \alpha\in Q_1}$  corresponding to $Q'$ is defined by
    \[
    M_i=\begin{cases}
        \mathbf{k}&\text{$i\in Q'_0$};\\
        0 &\text{otherwise}.
    \end{cases}\  \text{and}\  M_\alpha=\begin{cases}
        1_{\mathbf{k}} &\text{$\alpha\in Q_1'$};\\
        0&\text{otherwise}.
    \end{cases}
    \]
\end{theorem}

\subsection{Yoneda's Extension} Denote by $\mathbf{0}$ the zero object of $\opname{rep}(Q,\mathbb{F}_1)$.
An exact sequence of length $n+2$ of $\opname{rep}(Q,\mathbb{F}_1)$ is a sequence of morphism
\[
\xymatrix{L\ \ar[r]^{\alpha_0}&M_1\ar[r]^{\alpha_1}&M_2\ar[r]^{\alpha_2}&\cdots\ar[r]^{\alpha_{n-1}}&M_n\ar[r]^{\alpha_n}&N}
\]
 such that 
 \begin{itemize}
     \item $\ker \alpha_{i+1}=\opname{im} \alpha_i$ for $0\leq i\leq n-1$;
     \item $\alpha_0$ is a monomorphism, equivalently, $\ker\alpha_0=\mathbf{0}$;
     \item $\alpha_n$ is an epimorphism, equivalently, $\opname{im} \alpha_n=N$.
 \end{itemize}  As usual, we will use $\rightarrowtail$ (resp. $\twoheadrightarrow$) to indicate a morphism is a monomorphism  (resp. an epimorphism).

Let 
\[\xymatrix{\epsilon: &L\ \ar@{>->}[r]^{\alpha_0}&M_1\ar[r]^{\alpha_1}&M_2\ar[r]^{\alpha_2}&\cdots\ar[r]^{\alpha_{n-1}}&M_n\ar@{->>}[r]^{\alpha_n}&N}
\]
and 
\[\xymatrix{\epsilon': &U\ \ar@{>->}[r]^{\beta_0}&N_1\ar[r]^{\beta_1}&N_2\ar[r]^{\beta_2}&\cdots\ar[r]^{\beta_{n-1}}&N_n\ar@{->>}[r]^{\beta_n}&V}
\]
be  exact sequences of length $n+2$ of representations of $Q$ over $\mathbb{F}_1$. 
Denote by $\epsilon\oplus \epsilon'$ the following exact sequence 
\[
\xymatrix{L\oplus U\ \ar@{>->}[r]^-{\tiny \begin{pmatrix}\alpha_0&0\\0&\beta_0\end{pmatrix}}& M_1\oplus N_1\ar[r]&\cdots\ar[r]&M_n\oplus N_n\ar@{->>}[r]^-{\tiny\begin{pmatrix}\alpha_n&0\\0&\beta_n\end{pmatrix}}&N\oplus V,}
\]
which is called the direct sum of $\epsilon$ and $\epsilon'$.

Let $\mathbb{E}^n(N,L)$ be the set of exact sequences of length $n+2$ which is starting at $L$ and ending at $N$. An element in $\mathbb{E}^n(N,L)$ is called an $n$-extension of $N$ by L. We say that an element in $\mathbb{E}^n(N,L)$ is non-zero if at least one of its $n+2$ components is non-zero.
\begin{definition}
A non-zero exact sequence $\epsilon\in \mathbb{E}^n(N,L)$ is called primitive if $\epsilon$ has no non-trivial direct sum decomposition, that is, if $\epsilon= \epsilon_1\oplus \epsilon_2$ for exact sequences $\epsilon_1$ and $\epsilon_2$ of length $n+2$, then $\epsilon_1= 0$ or $\epsilon_2= 0$.
\end{definition}
\begin{remark}
Let $L,N\in \opname{rep}(Q,\mathbb{F}_1)$. Every non-zero exact sequence in $\mathbb{E}^n(N,L)$ is a finite direct sum of primitive exact sequences.
\end{remark}
Recall that the $n$-extensions $\epsilon$ and $\epsilon'$ of $N$ by $L$ satisfy the relation $\epsilon\leadsto\epsilon'$ if there is a commutative diagram
\[
\xymatrix{\epsilon: &L\ \ar@{=}[d]\ar@{>->}[r]&E_n\ar[d]^{f_n}\ar[r]&\cdots\ar[r]&E_1\ar[d]^{f_1}\ar@{->>}[r] &N\ar@{=}[d]\\
\epsilon':&L\ \ar@{>->}[r]&E_n'\ar[r]&\cdots\ar[r]&E_1'\ar@{->>}[r] &N.
}
\]
 If moreover $f_1,\dots, f_n$ are isomorphisms, we say that $\epsilon$ is isomorphic to $\epsilon'$, and denote it by $\epsilon\cong \epsilon'$.
The relation $\leadsto$ is not symmetric, but it generates an equivalence relation on $\mathbb{E}^n(N,L)$. We denote by $[\epsilon]$ the equivalence class of the $n$-extension $\epsilon$ and 
denote by $\Ext^n(N,L)$ the set of all equivalence classes of $n$-extensions of $N$ by $L$. The set $\Ext^n(N,L)$ is a pointed set with
\[0=[\xymatrix{L\ar@{=}[r]&L\ar[r]&\mathbf{0}\ar[r]&\cdots\ar[r]&\mathbf{0}\ar[r]&N\ar@{=}[r]&N}]
\] for $n\geq 2$ and 
\[\xymatrix{0=[L\ \ar@{>->}[r]^-{\tiny\begin{pmatrix}1_L&0\end{pmatrix}}&L\oplus N\ar@{->>}[r]^-{\tiny\begin{pmatrix}0\\ 1_N\end{pmatrix}}&N]}
\]
for $n=1$. By convention, we denote by $\Hom(L,M):=\Ext^0(L,M)$.
The following is obvious.
\begin{lemma}\label{l:equiv=0}
\begin{itemize}
    \item[(1)] For any exact sequences $\epsilon$ and $\delta$ of length $n+2$. If $[\epsilon]=0$ and $[\delta]=0$, then $[\epsilon\oplus \delta]=0$.
    \item[(2)] For any $n\geq 0$ and $L\in \opname{rep}(Q,\mathbb{F}_1)$, \[\Ext^n(\mathbf{0},L)=0=\Ext^n(L,\mathbf{0}).\] If $\Ext^n(L,-)=0$, then $\Ext^{n+1}(L,-)=0$.
\end{itemize}

\end{lemma}

\begin{definition}
   The category $\opname{rep}(Q,\mathbb{F}_1)$ is of global dimension $n$ if $\Ext^n(-,-)\neq 0$ and $\Ext^{n+1}(-,-)=0$. In this case, we denote it by $gl.dim \opname{rep}(Q,\mathbb{F}_1)=n$. It is hereditary if $gl.dim \opname{rep}(Q,\mathbb{F}_1)\leq 1$. 
\end{definition}

\begin{definition}
    A representation $P\in \opname{rep}(Q,\mathbb{F}_1)$ is projective if for any epimorphism $f:M\twoheadrightarrow N$  and any morphism $g:P\to N$ of $\opname{rep}(Q,\mathbb{F}_1)$, there is a morphism $h:P\to M$ such that $g=f\circ h$.
\end{definition}

The following is clear.
\begin{lemma}
Let $P\in \opname{rep}(Q,\mathbb{F}_1)$ be a projective representation.  Then $\Ext^1(P,-)=0$ and each direct summand of $P$ is projective.

\end{lemma}
\begin{remark}
In general, the condition $\Ext^1(P,-)=0$ does not imply that $P$ is projective and the direct sum of two projective representations is not projective, see example in Section \ref{s:non-linear}. 
\end{remark}

\section{Quivers of type $\mathbb{A}_n$ with linear orientations}\label{s:linear-quiver}
Throughout this section, let $\vec{\mathbb{A}}_n$ be the quiver of type $\mathbb{A}_n$ with a linear orientation, that is,
\[
\xymatrix{n\ar[r]^-{\alpha_{n-1}}&n-1\ar[r]^-{\alpha_{n-2}}&n-2\ar[r]&\cdots\ar[r]&2\ar[r]^-{\alpha_{1}}&1.}
\]
In this section, we study homological properties of $\opname{rep}(\vec{\mathbb{A}}_n,\mathbb{F}_1)$.

\subsection{Basic structure of $\opname{rep}(\vec{\mathbb{A}}_n,\mathbb{F}_1)$}
For $1\leq k\leq l\leq n$, denote by $[k,l]$ the $\mathbb{F}_1$-representation $M=(M_i,M_{\alpha_i})$ of $Q$  such that
\[M_i=\begin{cases}\mathbf{k} &\text{if $k\leq i\leq l$};\\
0& \text{otherwise}.
\end{cases} \text{and}\ M_{\alpha_i}=\begin{cases}1_{\mathbf{k}} &\text{if $k\leq i\leq l-1$};\\
0&\text{otherwise}.
\end{cases}
\]
Namely, $[k,l]$ is the indecomposable representation associated with the connected subquiver of $Q$ consisting of vertices $k,\dots, l$. By Theorem \ref{t:class-tree}, the isomorphism classes of indecomposable representations of $Q$ are precisely $\{[k,l]~|~1\leq k\leq l\leq n\}$. For each indecomposable representation $M\in \opname{rep}(\vec{\mathbb{A}}_n,\mathbb{F}_1)$, there is a unique simple representation $S_1$ (resp. $S_2$) such that $\Hom(S_1,M)\neq 0$ (resp. $\Hom(M,S_2)\neq 0$). In this case, we call $S_1$ and $S_2$ the {\it socle} and {\it top} of $M$ respectively.

We also need the following result on homomorphisms between indecomposable representations, which can be easily deduced from definition.
\begin{lemma}
Let $1\leq i\leq j\leq n$ and $1\leq k\leq l\leq n$.
We have
\begin{enumerate}
\item The socle of $[k,l]$ is the simple representation $[k,k]$, while, the top of $[k,l]$ is the simple representation  $[l,l]$.
    \item $\dim \Hom([i,j],[k,l])=\begin{cases}1& \text{if $i\leq k\leq j\leq l$;}\\ 0 &\text{otherwise.}
\end{cases}$
 \item Let $f\in \Hom([i,j],[k,l])$ be the unique non-zero homomorphism.
 \begin{itemize}
     \item[(a)] $f$ is an epimorphism if and only $i\leq k\leq j=l$;
     \item[(b)] $f$ is a monomorphism if and only if $i=k\leq j\leq l$;
     \item[(c)] $[k,k]$ is the socle of $\opname{im} f$.
 \end{itemize}
\end{enumerate} 
\end{lemma}

The following``splitting lemma" is an elementary observation, which plays a key role in the next subsection. 
\begin{lemma}\label{l:decom-epi}
Let $f:M\to N$ be an epimorphism. Let $N_1,N_2$ be subrepresentations of $N$ such that $N=N_1\oplus N_2$ and  $N_1\cong [k,l]$ for some $1\leq k\leq l\leq n$. Then there exist subrepresentations $M_1,M_2$ of $M$ such that 
\begin{enumerate}
\item $M=M_1\oplus M_2$;
\item $M_1\cong [i_0,l]$ for some $i_0\leq k$;
\item $f$ can be rewritten as $M_1\oplus M_2\xrightarrow{\tiny\begin{pmatrix}f_1&0\\ 0&f_2\end{pmatrix}}N_1\oplus N_2$, where $f_1$ is the unique non-zero morphism from $[i_0,l]$ to $[k,l]$ and $f_2$ is an epimorphism.
\end{enumerate}
\end{lemma}
\begin{proof}
    Denote by $M=(M_i, \phi_{\alpha_j})$, $N_1=(N_{1i},\varphi_{\alpha_j})$ and $N_2=(N_{2i},\varphi_{\alpha_j})$.
     Let $M_{1i}$ be the subspace of $M_i$ defined as follows
     \[M_{1i}=\begin{cases}0& \text{for $i>l$};\\
     f_l^t(N_{1l}) &\text{for $i=l$};\\
     \phi_{\alpha_i}(M_{1i+1}) &\text{for $i<l$}.
     \end{cases}
     \]
     It is straightforward to check that $M_1:=(M_{1i}, \phi_{\alpha_{j-1}}|_{M_{1j}})$ is a subrepresentation of $M$. Moreover, $M_2:=M\backslash M_1$ is also a subrepresentation of $M$  and $M=M_1\oplus M_2$. Since $f: M\to N$ is an epimorphism of representations and $N_1\cong [k,l]$, we conclude that $M_1\cong [i_0,l]$ for some $i_0\leq k$.
     Finally, by setting $f_1=f|_{M_1}$ and $f_2=f|_{M_2}$, we obtain the desired decomposition of $f$ by noting that $f$ is an epimorphism.
\end{proof}
\begin{remark}
A dual version of Lemma \ref{l:decom-epi} for monomorphisms also holds.
\end{remark}
\begin{corollary}\label{c:decom-seq}
Let $N=\bigoplus_{i=1}^tN_i$ with indecomposable direct summands $N_i$, $1\leq i\leq t$. Every exact sequence $\epsilon$ ending at $N$ of length $n+2$ has a decomposition
\[
\epsilon=\epsilon_1\oplus \cdots\oplus \epsilon_t,
\]
where $\epsilon_i$ is an exact sequence  of length $n+2$ ending at $N_i$.
\end{corollary}

\subsection{Global dimension}

In this section, we prove our main result on the global dimension of $\opname{rep}(\vec{\mathbb{A}}_n,\mathbb{F}_1)$.

\begin{lemma}\label{l:vanishi-E-3}
    $\Ext^3(-,-)=0$ for $\opname{rep}(\vec{\mathbb{A}}_n,\mathbb{F}_1)$.
\end{lemma}
\begin{proof}
 According to Corollary \ref{c:decom-seq} and Lemma \ref{l:equiv=0}, it suffices to show $\Ext^3(N,-)=0$ for any indecomposable representation $N$. Assume that $N=[k,l]$.
Let
\[\xymatrix{\epsilon: &L\ \ar@{>->}[r]^{f_0}&M_1\ar[r]^{f_1}&M_2\ar[r]^{f_2}&M_3\ar@{->>}[r]^{f_3}&N
}
\]
be an arbitrary exact sequence in $\mathbb{E}^3(N,L)$. By Lemma \ref{l:equiv=0} and Corollary \ref{c:decom-seq} again, we may assume that $\epsilon$ is primitive.
According to Lemma \ref{l:decom-epi}, $f_3$ admits a decomposition as 
\[\xymatrix{[i,l]\oplus\overline{M}_3\ar@{->>}[rr]^-{\tiny\begin{pmatrix}f_{31}&0\end{pmatrix}}&&[k,l],}\]
where $i\leq k$ and $f_{31}$ is the unique non-zero homomorphism from $[i,l]$ to $[k,l]$.
If $\overline{M}_3\neq 0$, then by applying Lemma \ref{l:decom-epi}, we will obtain a non-trivial decomposition of $\epsilon$, a contradiction. In particular, $M_3\cong [i,l]$.
Now we have $\ker f_3\cong [i,k-1]$\footnote{Here we use the convention that  $[s,t]$ is the zero representation whenever $t<s$.}. Similar to the preceding discussion, we conclude that $M_2\cong [i_1,k-1]$ for some $i_1\leq i$.  Now we have the following commutative diagram
\[\xymatrix{\delta:&L\ar@{=}[d]\ar@{=}[r]&L\ar[d]^{f_0}\ar[r]^-0 &[i_1,k-1]\ \ar@{=}[d]\ar@{>->}[r]&[i_1,l]\ar@{->>}[d]\ar@{->>}[r]^{\pi}&[k,l]\ar@{=}[d]\\
\epsilon:&
L\ \ar@{>->}[r]^{f_0}&M_1\ar[r]^-{f_1}&[i_1,k-1]\ar[r]^{f_2}&[i,l]\ar@{->>}[r]^{f_3}&[k,l].
}
\]
Consequently, $[\epsilon]=[\delta]$ in $\Ext^3(N,L)$. On the other hand, we also have the following commutative diagram
\[
\xymatrix{\delta:&L\ar@{=}[d]\ar@{=}[r]&L\ar@{=}[d]\ar[rr]^0 &&[i_1,k-1]\ \ar[d]^0\ar@{>->}[r]&[i_1,l]\ar@{->>}[d]^{\pi}\ar@{->>}[r]^{\pi}&[k,l]\ar@{=}[d]\\
\gamma:&L\ar@{=}[r]&L\ar[rr]^0 &&0\ar[r]^0&[k,l]\ar@{=}[r]&[k,l].
}
\]
We conclude that $[\epsilon]=0$ in $\Ext^3(N,L)$.
\end{proof}
\begin{lemma}\label{l:vanish-2}
    Assume that $n\leq 2$. Then $\Ext^2(-,-)=0$ for $\opname{rep}(\vec{\mathbb{A}}_n,\mathbb{F}_1)$.
\end{lemma}
\begin{proof}
    This is clear true for $n=1$.  Let us consider the case $n=2$. According to Corollary \ref{c:decom-seq}, it suffices to show $\Ext^2(L,-)=0$ for any indecomposable representation $L$. This is also clear for $L=[1,1]$ or $[1,2]$. Now assume that $L=[2,2]$. 
    Let \[\xymatrix{\delta: U\ \ar@{>->}[r]&M\ar[r]&N\ar@{->>}[r]&[2,2]}\] be an arbitrary exact sequence. By applying Lemma \ref{l:decom-epi}, there is an $i$ such that $N=[i,2]\oplus \overline{N}$. If $i=2$, then $\delta$ is the direct sum of \[\xymatrix{U\ \ar@{>->}[r]&M\ar[r]&\overline{N}\ar@{->>}[r]& \mathbf{0}}\] and \[\xymatrix{\mathbf{0}\ \ar@{>->}[r]&\mathbf{0}\ar[r]&[2,2]\ar@{=}[r]&[2,2]}.\]
     If $i=1$, then $\delta$ is a direct sum of \[\xymatrix{U\ \ar@{>->}[r]&\overline{M}\ar[r]&\overline{N}\ar@{->>}[r]& \mathbf{0}}\] and \[\xymatrix{\mathbf{0}\ \ar@{>->}[r]&[1,1]\ar[r]&[1,2]\ar@{->>}[r]&[2,2]}.\]  In either case, we have $[\delta]=0$.
\end{proof}
\begin{remark}
    In fact, one can show that $\Ext^1(-,-)=0$ for $n=1$. On the other hand, for $\epsilon: \xymatrix{[1,1]\ \ar@{>->}[r]&[1,2]\ar@{->>}[r]&[2,2]}$, we have $[\epsilon]\neq 0$. Hence $\Ext^1(-,-)\neq 0$ for $n=2$.
\end{remark}
\begin{lemma}\label{l:E-3-split}
    Assume that $n\geq 3$. Let \[\delta: \xymatrix{[1,1]\ \ar@{>->}[r]&M\ar[r] &N\ar@{->>}[r]&[3,3]}\] be an exact sequence in $\mathbb{E}^2([3,3],[1,1])$ such that there is a commutative diagram
    \[\xymatrix{\delta: &[1,1]\ \ar@{=}[d]\ar@{>->}[r]&M\ar[d]\ar[r] &N\ar[d]\ar@{->>}[r]&[3,3]\ar@{=}[d]\\
    \epsilon:& [1,1]\ \ar@{>->}[r]&[1,2]\ar[r]&[2,3]\ar@{->>}[r]&[3,3]}
    \]
    or 
    \[\xymatrix{\delta:&[1,1]\ \ar@{=}[d]\ar@{>->}[r]&M\ar[r] &N\ar@{->>}[r]&[3,3]\ar@{=}[d]\\
    \epsilon:&[1,1]\ \ar@{>->}[r]&[1,2]\ar[u]\ar[r]&[2,3]\ar[u]\ar@{->>}[r]&[3,3].}
    \]
    Then $\delta\cong\epsilon\oplus \gamma$, where $\gamma: \xymatrix{\mathbf{0}\ \ar@{>->}[r]&U\ar@{=}[r] &U\ar@{->>}[r]&\mathbf{0}}$ for some $U$.
\end{lemma}
\begin{proof}
    Let us prove the first case, the second case can be proved in a similar way. Assume that we have a commutative diagram 
    \[\xymatrix{\delta: &[1,1]\ \ar@{=}[d]\ar@{>->}[r]&M\ar[d]\ar[r] &N\ar[d]\ar@{->>}[r]&[3,3]\ar@{=}[d]\\
    \epsilon:& [1,1]\ \ar@{>->}[r]&[1,2]\ar[r]&[2,3]\ar@{->>}[r]&[3,3].}
    \]
    Applying Lemma \ref{l:decom-epi}, we conclude that there is an $1\leq i\leq 2$ such that we can rewrite the commutative diagram as
    \[\xymatrix{\delta: &[1,1]\ \ar@{=}[d]\ar@{>->}[r]&M\ar[d]\ar[r] &[i,3]\oplus \overline{N}\ar[d]^{\tiny\begin{pmatrix}
        f&0
    \end{pmatrix}}\ar@{->>}[r]^-{\tiny\begin{pmatrix}
        \pi&0
    \end{pmatrix}}&[3,3]\ar@{=}[d]\\
    \epsilon:& [1,1]\ \ar@{>->}[r]&[1,2]\ar[r]&[2,3]\ar@{->>}[r]&[3,3],}
    \]
    where $\pi:[i,3]\to [3,3]$ and $f:[i,3]\to [2,3]$ are the unique non-zero morphisms respectively.

    We claim that $i=2$ and hence $f$ is the identity map. Otherwise, $i=1$. In this case, $\ker (\pi\ 0)=[1,2]\oplus \overline{N}$ and we have an exact sequence 
    \[\xymatrix{[1,1]\ \ar@{>->}[r]&M\ar@{->>}[r]&[1,2]\oplus \overline{N}.}
    \]
    By Lemma \ref{l:decom-epi} again, we can rewrite the above short exact sequence as
    \[\xymatrix{[1,1]\ \ar@{>->}[rr]^-{\tiny\begin{pmatrix}
        0\\ l
    \end{pmatrix}}&&[1,2]\oplus \overline{M}\ar@{->>}[rr]^-{\tiny\begin{pmatrix}
        1_{[1,2]} &0\\0&h
    \end{pmatrix}}&&[1,2]\oplus \overline{N}.}
    \]
    Consequently, the original commutative diagram can be rewritten as 
    \[\xymatrix{\delta: &[1,1]\ \ar@{=}[d]\ar@{>->}[r]^-{\tiny\begin{pmatrix}
        0\\ l
    \end{pmatrix}}&[1,2]\oplus \overline{M}\ar[d]^{\tiny\begin{pmatrix}
        0&\beta
    \end{pmatrix}}\ar[r]^-{\tiny\begin{pmatrix}
        g &0\\0&h
    \end{pmatrix}} &[1,3]\oplus \overline{N}\ar[d]^{\tiny\begin{pmatrix}
        f&0
    \end{pmatrix}}\ar@{->>}[r]^-{\tiny\begin{pmatrix}
        \pi&0
    \end{pmatrix}}&[3,3]\ar@{=}[d]\\
    \epsilon:& [1,1]\ \ar@{>->}[r]&[1,2]\ar[r]&[2,3]\ar@{->>}[r]&[3,3],}
    \]
    where $g:[1,2]\to [1,3]$ is the unique non-zero morphism. However, the middle square is not commutative, a contradiction. Hence $i=2$ and $\ker (\pi\ 0)=[2,2]\oplus \overline{N}$. By applying Lemma \ref{l:decom-epi}, we conclude that there is a $1\leq j\leq 2$ such that $M=[j,2]\oplus \overline{M}$. A similar discussion shows that $j=1$ and the commutative diagram can be rewritten as
    \[\xymatrix{\delta: &[1,1]\ \ar@{=}[d]\ar@{>->}[r]^-{\tiny\begin{pmatrix}
        l\\ 0
    \end{pmatrix}}&[1,2]\oplus \overline{M}\ar[d]^{\tiny\begin{pmatrix}
        1_{[1,2]}&0
    \end{pmatrix}}\ar[r]^-{\tiny\begin{pmatrix}
        g &0\\0&u
    \end{pmatrix}} &[2,3]\oplus \overline{N}\ar[d]^{\tiny\begin{pmatrix}
        1_{[2,3]}&0
    \end{pmatrix}}\ar@{->>}[r]^-{\tiny\begin{pmatrix}
        \pi&0
    \end{pmatrix}}&[3,3]\ar@{=}[d]\\
    \epsilon:& [1,1]\ \ar@{>->}[r]&[1,2]\ar[r]&[2,3]\ar@{->>}[r]&[3,3].}
    \]
     Moreover, $u:\overline{M}\to \overline{N}$ is an isomorphism. We conclude that $\delta\cong\epsilon\oplus \gamma$, where $\gamma:\xymatrix{\mathbf{0}\ \ar@{>->}[r]&\overline{N}\ar@{=}[r] &\overline{N}\ar@{->>}[r]&\mathbf{0}}$.
\end{proof}

\begin{corollary}\label{c:non-vanish3}
    Assume that $n\geq 3$. Then $\Ext^2(-,-)\neq 0$ for $\opname{rep}(\vec{\mathbb{A}}_n,\mathbb{F}_1)$.
\end{corollary}
\begin{proof}
Let \[\epsilon: \xymatrix{[1,1]\ \ar@{>->}[r]&[1,2]\ar[r]&[2,3]\ar@{->>}[r]&[3,3]}\]
be the exact sequence in $\mathbb{E}^{2}([3,3],[1,1])$.
    According to Lemma \ref{l:E-3-split} and the definition of the equivalence relation generated by $\leadsto$, we have $[\epsilon]\neq 0$ in $\Ext^2([3,3],[1,1])$.
\end{proof}
\begin{theorem}\label{t:main-result}
The category $\opname{rep}(\vec{\mathbb{A}}_n,\mathbb{F}_1)$ is hereditary if and only if $n\leq 2$. Moreover, $\opname{rep}(\vec{\mathbb{A}}_n,\mathbb{F}_1)$ is of global dimension $2$ whenever $n\geq 3$.
\end{theorem}
\begin{proof}
    The first statement follows from Lemma \ref{l:vanish-2} and Corollary \ref{c:non-vanish3}. The second statement follows from Lemma \ref{l:vanishi-E-3} and Corollary \ref{c:non-vanish3}.
\end{proof}

\subsection{Projective representations}

\begin{lemma}\label{l:projective}
    Let $P=(P_i,P_\alpha)\in \opname{rep}(\vec{\mathbb{A}}_n,\mathbb{F}_1)$ such that $P_\alpha$ is an injection for each arrow $\alpha$. Then $P$ is a projective representation.
\end{lemma}
\begin{proof}
    Let $M=(M_i,M_\alpha)$ and $N=(N_i,N_\alpha)$ be $\mathbb{F}_1$-representations of $\vec{\mathbb{A}}_n$ and $f=(f_i):M\twoheadrightarrow N$ be an epimorphism. Let $g=(g_i):P\to N$ be any morphism. We have to show that there is a morphism $h=(h_i):P\to M$ such that $g=f\circ h$. We are going to construct $h_i$ inductively. For $n$, we define
    \[h_n(p_n)=\begin{cases}0&\text{if $p_n\in \ker g_n$;}\\
    f_n^t(g_n(p_n)) & \text{if $p_n\not\in \ker g_n$.}
    \end{cases}
    \]
    It is obvious that $h_n$ is an $\mathbb{F}_1$-linear map and $g_n=f_n\circ h_n$.

    Now suppose that $h_j$ is defined such that $g_j=f_j\circ h_j$. We define 
    \[
    h_{j-1}(p_{j-1})=\begin{cases}
        M_{\alpha_{j-1}}(h_j(P_{\alpha_{j-1}}^t(p_{j-1}))) &\text{if $p_{j-1}\in \opname{im} P_{\alpha_{j-1}}$;}\\
0&\text{if $p_{j-1}\not\in \opname{im} P_{\alpha_{j-1}}$, $p_{j-1}\in \ker g_{j-1}$;}\\
f_{j-1}^t(g_{j-1}(p_{j-1}))& \text{if $p_{j-1}\not\in \opname{im} P_{\alpha_{j-1}}$, $p_{j-1}\not\in \ker g_{j-1}$;}
    \end{cases}
    \]
    In order to show that $h_{j-1}$ is an $\mathbb{F}_1$-linear map, it suffices to show that \[M_{\alpha_{j-1}}\circ h_j\circ P_{\alpha_{j-1}}^t(x_{j-1})\neq f_{j-1}^t(g_{j-1}(y_{j-1}))\] for $x_{j-1}\in \opname{im} P_{\alpha_{j-1}}$ and $y_{j-1}\not\in \opname{im} P_{\alpha_{j-1}}$, $y_{j-1}\not\in \ker g_{j-1}$. Otherwise,
    \[
    f_{j-1}(M_{\alpha_{j-1}}\circ h_j\circ P_{\alpha_{j-1}}^t(x_{j-1}))=f_{j-1}(f_{j-1}^t(g_{j-1}(y_{j-1})))=g_{j-1}(y_{j-1})\neq 0.
    \]
    On the other hand, we have 
    \begin{eqnarray*} f_{j-1}(M_{\alpha_{j-1}}\circ h_j\circ P_{\alpha_{j-1}}^t(x_{j-1}))&=&N_{\alpha_{j-1}}\circ f_j\circ h_j\circ P_{\alpha_{j-1}}^t(x_{j-1})\\
    &=&N_{\alpha_{j-1}}\circ g_j\circ P_{\alpha_{j-1}}^t(x_{j-1})\\
    &=&g_{j-1}\circ P_{\alpha_{j-1}}\circ P_{\alpha_{j-1}}^t(x_{j-1})\\
    &=&g_{j-1}(x_{j-1}),
    \end{eqnarray*}
    which contradicts to that $g_{j-1}$ is an $\mathbb{F}_1$-linear map. It is straightforward to check that $g_{j-1}=f_{j-1}\circ h_{j-1}$. Finally, for any $p_j\in P_j$, we have 
    \[
    h_{j-1}(P_{\alpha_{j-1}}(p_j))=M_{\alpha_{j-1}}\circ h_j\circ P_{\alpha_{j-1}}^t(P_{\alpha_{j-1}}(p_j))=M_{\alpha_{j-1}}\circ h_j(p_j),
    \]
    where the last equality follows from the fact that $P_{\alpha_{j-1}}$ is an injection.
\end{proof}
\begin{proposition}
    A representation $M=(M_i,M_\alpha)$ is projective if and only if $M_\alpha$ is an injection for each arrow $\alpha$ of $\vec{\mathbb{A}}_n$.
\end{proposition}
\begin{proof}
    According to Lemma \ref{l:projective}, it remains to show that if $M=(M_i,M_\alpha)$ is projective, then $M_\alpha$ is an injection for  each arrow $\alpha$. Since the direct summand of a projective representation is a projective, it suffices to assume that $M$ is indecomposable. Namely, we may assume that $M=[i,l]$ for some $1\leq i\leq l\leq n$. Clearly, we have an epimorphism $f:[1,l]\twoheadrightarrow [i,l]$. Now $[i,l]$ is projective implies that $[i,l]$ is a direct summand of $[1,l]$. Hence $[i,l]\cong [1,l]$ and we are done.
\end{proof}
\begin{corollary}\label{c:proj-direct-sum}
Direct sum of projective representations in $\opname{rep}(\vec{\mathbb{A}}_n,\mathbb{F}_1)$ is projective. 
\end{corollary}
Recall that we say $\opname{rep}(Q,\mathbb{F}_1)$ has {\it enough projectives} if for any $M\in \opname{rep}(Q,\mathbb{F}_1)$, there is a projective representation $P_M$ and an epimorphism $f:P_M\twoheadrightarrow M$. The following is a direct consequence of Theorem \ref{t:class-tree}, Lemma \ref{l:projective} and Corollary \ref{c:proj-direct-sum}.
\begin{corollary}
The category $\opname{rep}(\vec{\mathbb{A}}_n,\mathbb{F}_1)$ has enough projectives.
\end{corollary}

\section{(Counter)-examples}\label{s:exam}
\subsection{Euler form}
Let $\vec{\mathbb{A}}_n$ be the linear quiver as Section \ref{s:linear-quiver}. According to Theorem \ref{t:main-result} and Corollary \ref{c:decom-seq}, one can show that for any $L,N\in \opname{rep}(\vec{\mathbb{A}}_n,\mathbb{F}_1)$, $\Ext^i(L,N)$ is a finite point set for each $i\geq 0$. Consequently,  the Euler form $\langle-, -\rangle:\opname{rep}(\vec{\mathbb{A}}_n,\mathbb{F}_1)\times \opname{rep}(\vec{\mathbb{A}}_n,\mathbb{F}_1)\to \mathbb{Z}$ is well-defined, where
\[
\langle L,N\rangle:=\sum_{i=0}^\infty(-1)^i\dim\Ext^i(L,N),
\]
for $L,N\in \opname{rep}(\vec{\mathbb{A}}_n,\mathbb{F}_1)$.

Now assume moreover that $n=2$.  According to Theorem \ref{t:class-tree}, there are precisely $3$ indecomposable $\mathbb{F}_1$-representations:
\[
\xymatrix{S_1: 0\ar[r]&\mathbf{k}}, \xymatrix{S_2: \mathbf{k}\ar[r]&0}, \xymatrix{P_2: \mathbf{k}\ar[r]^-{1_{\mathbf{k}}}&\mathbf{k}}.
\]
According to Theorem \ref{t:main-result}, $\Ext^2(-,-)=0$ for $\opname{rep}(\vec{\mathbb{A}}_2,\mathbb{F}_1)$. Thus for any $L,M\in \opname{rep}(\vec{\mathbb{A}}_2,\mathbb{F}_1)$, \[\langle L,M\rangle=\dim\Hom(L,M)-\dim\Ext^1(L,M).\]
Note that the Grothendieck group $\opname{G}_0(\opname{rep}(\vec{\mathbb{A}}_2,\mathbb{F}_1))$ of $\opname{rep}(\vec{\mathbb{A}}_2,\mathbb{F}_1)$ identifies with $\mathbb{Z}^2$ via $\dimv$. We are going to show that $\langle P_2\oplus P_2,P_2\rangle \neq \langle S_1\oplus S_2\oplus P_2, S_1\oplus S_2\rangle$, which implies that $\langle-,-\rangle$ does not descend to $\opname{G}_0(\opname{rep}(\vec{\mathbb{A}}_2,\mathbb{F}_1))$.

Since $P_2\oplus P_2$ is projective by Lemma \ref{l:projective}, we have \[\langle P_2\oplus P_2,P_2\rangle=\dim \Hom(P_2\oplus P_2, P_2)=2.\]
On the other hand, it is easy to see that $\dim\Hom(S_1\oplus S_2\oplus P_2, S_1\oplus S_2)=5$.
Since $S_1$ and $P_2$ are projective representations, every exact sequence in $\mathbb{E}^1(S_1\oplus S_2\oplus P_2, S_1\oplus S_2)$ is a direct sum of an exact sequence in $\mathbb{E}^1(S_2, S_1\oplus S_2)$ with $\xymatrix{\mathbf{0}\ \ar@{>->}[r]&S_1\oplus P_2\ar@{=}[r]&S_2\oplus P_2}$.  By Lemma \ref{l:decom-epi}, we conclude that there are exactly two non-isomorphic $1$-extension of $S_1\oplus S_2\oplus P_2$ by $S_1\oplus S_2$ up to isomorphism.
It is straightforward to see that these two elements are not equivalent and yield distinct elements in $\Ext^1(S_1\oplus S_2\oplus P_2, S_1\oplus S_2)$. Thus $\dim \Ext^1(S_1\oplus S_2\oplus P_2, S_1\oplus S_2)=1$ and 
\[
\langle S_1\oplus S_2\oplus P_2,S_1\oplus S_2\rangle=4.
\]

\subsection{Projectives in non-linear quiver}\label{s:non-linear}
Let $Q$ the following quiver $\xymatrix{3\ar[r]&2&1\ar[l]}$. According to Theorem \ref{t:class-tree}, $Q$ has precisely $6$ indecomposable $\mathbb{F}_1$-representations:
\[\xymatrix{S_3: \mathbf{k}\ar[r]&0&0\ar[l]}, \xymatrix{S_2: 0\ar[r]&\mathbf{k}&0\ar[l]}, \xymatrix{S_1: 0\ar[r]&0&\mathbf{k}\ar[l]}
\]
\[
\xymatrix{P_3: \mathbf{k}\ar[r]^-{1_{\mathbf{k}}}&\mathbf{k}&0\ar[l]}, \xymatrix{P_1:0\ar[r]&\mathbf{k}&\mathbf{k}\ar[l]_-{1_{\mathbf{k}}}}, 
\]
\[\xymatrix{M:\mathbf{k}\ar[r]^-{1_{\mathbf{k}}}&\mathbf{k}&\mathbf{k}\ar[l]_-{1_{\mathbf{k}}}}.
\]
It is straightforward to check  the following:
\begin{itemize}
    \item[(1)] $S_2,P_1,P_3$ are projective representations.
    \item[(2)] Let $N$ be an indecomposable representation such that $\Hom(N,M)\neq 0$. Then $N\in \{S_2,P_1,P_3,M\}$ and $\dim \Hom(N,M)=1$.
    \item[(3)] $\dim\Hom(M,P_1)=0$.
    \item[(4)] $\dim\Hom(P_1,S_1)=1$ and the non-zero morphism $f$ is an epimorphism.
    \item[(5)] $\dim\Hom(M,S_1)=1$ and we denote the non-zero morphism by $g$.
\end{itemize}
It follows that there is no morphism $h\in \Hom(M,P_1)$ such that $f\circ h=g$. Hence $M$ is not projective. 

On the other hand, assume that $u:X\twoheadrightarrow M$ is an epimorphism. Then it is easy to show that $X\cong M\oplus \overline{X}$ for some representation $\overline{X}$ and $u$ can be rewritten as $M\oplus \overline{X}\xrightarrow{(1_M \ 0)} M$. In particular, $\Ext^1(M,-)=0$. Furthermore, it is easy to see that the projection $P_1\oplus P_3\to S_1\oplus S_3$  does not factor along the epimorphism $M\twoheadrightarrow S_1\oplus S_3$, which implies that $P_1\oplus P_3$ is not projective.
We conclude that 
\begin{itemize}
    \item $\opname{rep}(Q,\mathbb{F}_1)$ has no enough projective representations;
    \item $\Ext^1(L,-)=0$ is not equivalent to that $L$ is a projecitve representation;
    \item the direct sum of projectives is not projective in general.
\end{itemize}

\end{document}